\documentclass[aap,MSNbibl,dvips]{arximspdf}
\usepackage{mathbh}

\doi{10.1214/12-AAP904} 
\volume{23}
\issue{6}
\pubyear{2013}
\firstpage{2357}
\lastpage{2381}

\makeatletter
\newcommand{\rrvert}{\vert}
\newcommand{\llvert}{\vert}
\newcommand{\eqref}[1]{(\ref{#1})}
\newcommand{\norm}[1]{\Vert #1 \Vert}
\newcommand{\E}{\mathbb{E}}
\renewcommand{\P}{\mathbb{P}}
\newcommand{\R}{\mathbb{R}}
\newcommand{\PX}{\mathbf{P}}
\newcommand{\Ccal}{\mathcal{C}}
\newcommand{\Ecal}{\mathcal{E}}
\newcommand{\Jcal}{\mathcal{J}}
\newcommand{\Lcal}{\mathcal{L}}
\newcommand{\Mcal}{\mathcal{M}}
\newcommand{\Ncal}{\mathcal{N}}
\newcommand{\Scal}{\mathcal{S}}
\newcommand{\Wcal}{\mathcal{W}}
\newcommand{\indicator}[1]{\mathbh{1}_{\{#1\}}}

\newtheorem{theorem}{Theorem}[section]
\newtheorem{prop}[theorem]{Proposition}
\newtheorem{lemma}[theorem]{Lemma}
\newtheorem{corollary}[theorem]{Corollary}
\newproclaim{assumption}[theorem]{Assumption}
\newproclaim{notation}[theorem]{Notation}
\newproclaim{defn}{Definition}
\makeatother

\begin{document}
\begin{frontmatter}

\title{Scaling limits via excursion theory: Interplay between
Crump--Mode--Jagers branching processes and processor-sharing queues}
\runtitle{Scaling limits via excursion theory}

\begin{aug}
\author[A]{\fnms{Amaury} \snm{Lambert}\ead[label=e1]{amaury.lambert@upmc.fr}\thanksref{t0}},
\author[B]{\fnms{Florian} \snm{Simatos}\corref{}\ead[label=e2]{florian.simatos@cwi.nl}\thanksref{t1}}
\and
\author[C]{\fnms{Bert} \snm{Zwart}\ead[label=e3]{bert.zwart@cwi.nl}\thanksref{t1}}
\thankstext{t0}{Funded by project MANEGE 09-BLAN-0215 of ANR (French
national research agency).}
\thankstext{t1}{Sponsored by an NWO-VIDI grant.}
\runauthor{A. Lambert, F. Simatos and B. Zwart}
\affiliation{UPMC Univ Paris 06, CWI and EURANDOM, and CWI, EURANDOM, VU~University Amsterdam and GeorgiaTech}
\address[A]{A. Lambert\\
Laboratoire de Probabilit\'es et Mod\`eles Al\'eatoires\\
UMR 7599 CNRS and UPMC Univ Paris 06\\
Case courrier 188\\
4 Place Jussieu\\
F-75252 Paris Cedex 05\\
France\\
\printead{e1}}
\address[B]{F. Simatos\\
CWI\\
123 Science Park\\
1098 XG, Amsterdam\\
The Netherlands\\
\printead{e2}}
\address[C]{B. Zwart\\
CWI\\
123 Science Park\\
1098 XG, Amsterdam\\
The Netherlands\\
\printead{e3}}
\end{aug}

\received{\smonth{4} \syear{2011}}
\revised{\smonth{10} \syear{2012}}

%
\begin{abstract}
We study the convergence of the $M/G/1$ processor-sharing,
queue length process in the heavy traffic regime, in the finite
variance case. To do so, we combine results pertaining to L\'evy
processes, branching processes and queuing theory. These results yield
the convergence of long excursions of the queue length processes,
toward excursions obtained from those of some reflected Brownian motion
with drift, after taking the image of their local time process by the
Lamperti transformation. We also show, via excursion theoretic
arguments, that this entails the convergence of the entire processes to
some (other) reflected Brownian motion with drift. Along the way, we
prove various invariance principles for homogeneous, binary
Crump--Mode--Jagers processes.
In the last section we discuss potential implications of the state
space collapse property, well known in the queuing literature, to
branching processes.
\end{abstract}

%


\end{frontmatter}

\section{Introduction}\label{sec1}

The standard machinery to show weak convergence of stochastic processes
consists in proving tightness and characterizing accumulation points.
Probably the most common technique to characterize accumulation points
is to show that finite-dimensional distributions converge, but as Jacod
and Shiryaev~\cite{Jacod030} point out, this is ``very often [\ldots]
a very difficult (or simply impossible) task to accomplish.'' In the
present work, motivated by the processor-sharing~(PS) queue length
process, we develop new ideas to characterize limit points of a
sequence of regenerative processes. The basic\vadjust{\goodbreak} idea is to show that the
convergence of suitably conditioned excursions implies the convergence
of the full processes. Our starting point to control excursions is the
Lamperti transformation that links excursions of the PS queue to
Crump--Mode--Jagers (CMJ) branching processes. Further, control on CMJ
processes comes from a recent result of Lambert~\cite{Lambert100} that
relates them to L\'evy processes via local times.

\subsection*{The processor-sharing queue}

The PS queue is a single-server queue in which the server splits its
service capacity equally among all the users present. For instance, if
the server has a service capacity $c$ and if there are exactly $q \geq
1$ customers in the queue during the time interval $[t, t+h]$ (in
particular there is no arrival or departure), then the residual service
requirement of each customer is decreased by $ch/q$ during this time
interval while the total workload is decreased by~$ch$.

\subsection*{Crump--Mode--Jagers branching processes}

A CMJ process is a stochastic process counting the size of a population
where individuals give birth to independent copies of themselves. It is
defined through a pair $(V, \xi)$ of possibly dependent random
variables, where $V > 0$ is a real valued random variable and $\xi$ is
a point process on $(0,\infty)$ (in particular, its atoms have
integer-valued weights). Each individual $a$ of the branching process
is given an independent copy $(V_a, \xi_a)$ of $(V, \xi)$: if the
individual $a$ is born at time $B_a$, then at time $B_a \leq t \leq B_a
+ V_a$ she gives birth to $\xi_a(\{t-B_a\})$ i.i.d. copies of herself,
$V_a$ therefore being seen as her life length. A CMJ process is called
binary and homogeneous when $\xi$ is a Poisson process independent from
$V$. Lambert~\cite{Lambert100} has shown that a binary, homogeneous
CMJ process which is in addition critical or subcritical is the local
time process of a suitable spectrally positive L\'evy process.

\subsection*{The Lamperti transformation}

Connections between branching processes and queues have been known for
a long time. Kendall~\cite{Kendall510} is usually referred to as one
of the earliest publications in this area. Concerning the PS queue,
Kitayev and Yashkov~\cite{Kitayev790} have proved that a busy cycle of
the PS queue length process becomes a CMJ process after a suitable time
change. This time-change transformation is the same one that links
continuous-state space branching processes and L\'evy processes and is
called Lamperti transformation in the branching literature; see
Lamperti~\cite{Lamperti671} or Caballero et al.~\cite{Caballero090}
for a more recent treatment. This connection between the PS queue and
CMJ processes has been used to establish results on the stationary
behavior of the PS queue; see, for instance, Grishechkin~\cite
{Grishechkin920}. In this paper we make a deeper use of this
connection, since we exploit it to study the entire trajectories of the
processes.

The connections between CMJ processes, L\'evy processes and PS queues
lead to a natural proof of the weak convergence of CMJ processes. On
the one hand, we can prove\vadjust{\goodbreak} tightness of such processes by transferring,
via the Lamperti transformation, a result in queueing theory on the
departure process of queues with a symmetric service discipline. On the
other hand, exploiting the fact that subcritical, binary and
homogeneous CMJ processes are local time processes of suitable L\'evy
processes makes it possible to characterize accumulation points.

Using continuity properties of the Lamperti transformation, much in the
spirit of those established by Helland~\cite{Helland780}, and the
connection between CMJ processes and the PS queue, the convergence of
suitably renormalized CMJ processes implies that excursions of the
($M/G/1$) PS queue length processes converge. Thus convergence of
excursions of the PS queue length process comes quite naturally by
combining different results from queueing theory and the theory of L\'
evy and CMJ processes. Besides the original combination of these
various results, the main methodological contribution of the present
work is to show that from there, one can conclude that the whole PS
queue length processes converge.

\subsection*{Weak convergence of CMJ processes}

Binary, homogeneous CMJ processes considered in the present paper can
be seen as branching processes and also as local time processes. Since
Lamperti~\cite{Lamperti670}, we have a complete characterization of
the possible asymptotic behaviors of branching processes in discrete
time. Grimvall~\cite{Grimvall730} improved Lamperti's results by
proving tightness and hence weak convergence; see also Chapter~$9$ in
Ethier and Kurtz~\cite{Ethier860} for another proof of Grimvall's
results using time-change arguments. In the continuous-time setting,
the Markovian case has been studied by Helland~\cite{Helland780}.
Outside the Markovian case we are only aware of two papers by
Sagitov~\cite{Sagitov940,Sagitov950}, that establish convergence of
the finite-dimensional distributions for some particular CMJ that are
not homogeneous.

As for convergence of local time processes, there is a wealth of
literature studying the convergence of local time processes associated
to random walks converging to Brownian motion. One of the earliest
publication in this domain is Knight~\cite{Knight630}; see also
Borodin~\cite{Borodin810,Borodin860}, Perkins~\cite{Perkins820} and
references therein. The problem of finding sharp convergence rates has
been the focus of intense activity; see, for instance, the introduction
of Cs\"{o}rg\H{o} and R\'ev\'esz~\cite{Csorgo850} for references. On
the other hand, the question of the convergence of local time processes
associated to compound Poisson processes (which is another natural way
to approximate a Brownian motion) has been comparatively very little
studied. In this context, Khoshnevisan~\cite{Khoshnevisan930} derived
sharp convergence rates using embedding techniques requiring bounded
fourth moments. In that respect, some of the results of the present
paper seem to be new. Under a second moment assumption on lifetimes,
Theorem~\ref{thmconv-CMJ-excursion} shows the weak convergence of
suitably renormalized homogeneous, binary CMJ processes when started
from one individual and conditioned by their total offspring.
Theorem~\ref{thmconv-CMJ-star} states the weak convergence of these
CMJ processes when started\vadjust{\goodbreak} with large initial condition to the Feller
diffusion. In this last setting, the choice of the initial condition
turns out to be quite subtle, and is discussed in Section~\ref{secdiscussion}.

\subsection*{Heavy traffic of the PS queue}

The heavy traffic limit of the PS queue has been investigated by
Gromoll~\cite{Gromoll040}. Extending the state space collapse
framework developed by Bramson~\cite{Bramson980} and Williams~\cite
{Williams980}, he proved convergence of the measure-valued descriptors
of the $G/G/1$-PS queue, assuming that service requirements have
bounded fourth moments, toward some measure-valued diffusion process.
In the present paper we assume Poisson arrivals, that is, we study the
$M/G/1$-PS queue, but we relax the moment assumption and prove
convergence of the queue length process under a minimal second moment
assumption. More specifically, Theorem~\ref{thmconv-PS} shows the
convergence of these (suitably renormalized) queue length processes to
the regenerative process whose excursions are obtained from those of
some reflected Brownian motion with drift after taking the image of
their local time process by the Lamperti transformation. Theorem~\ref
{thmconv-PS-2} shows that this process actually is another reflected
Brownian motion with drift.

We also believe that our method can be used to study the case of
service requirements with infinite variance, where state space collapse
cannot be used since the workload and queue length processes have
different orders of magnitude. To the very least, although so far there
was no conjecture to this open problem, our method clearly suggests a
candidate for the heavy traffic limit of the PS queue length process in
the infinite variance case. Following the arguments in the previous
paragraphs, this limit should be a regenerative process whose
excursions away from 0 are obtained from those of some reflected,
spectrally positive L\'evy process by taking first their local time
process and then applying Lamperti transformation.

\subsection*{Organization of the paper}

In Section~\ref{secpreliminaries} we introduce general notation and
state preliminary results. In Section~\ref{secCMJ-PS-Levy}, we explain
the connections between CMJ processes, PS queues and L\'evy processes.
We also introduce the processes studied throughout the rest of the
paper. Section~\ref{secmain} is devoted to the proof of the main
result of the paper, Theorem~\ref{thmconv-PS}, which states the
convergence of the PS queue length process toward a process that we
define through its excursion measure. In Section~\ref
{secinitial-state} we extend this result by explicitly identifying the
limiting process as being a(nother) reflected Brownian motion with
drift and by considering a general initial condition. Finally, in
Section~\ref{secdiscussion} we make some comments about continuity
properties of local time processes and possible implications of the
state-space collapse property to branching processes.

\section{Notation and preliminary results} \label{secpreliminaries}

Let $D$, respectively, $D_+$, be the set of c\`adl\`ag functions from
$[0,\infty)$ to $\R$, respectively, to $[0,\infty)$. For $f \in D$ and
$m \geq0$ let $\norm{f}_m = \sup_{[0,m]} |f|$\vadjust{\goodbreak} and $\norm{f}_\infty=
\sup|f|$. We will endow $D$ with the topology of uniform convergence
on compact sets, that is, we will write $f_n \to f$ for functions $f_n,
f \in D$ if $\norm{f_n - f}_m \to0$ as $n \to\infty$ for every $m
\geq0$. The space $D$ is more naturally endowed with the Skorohod
$J_1$ topology (see, e.g., Billingsley~\cite{Billingsley990}),
but since the Skorohod topology relativized to the space of continuous
functions coincides with the topology of uniform convergence on compact
sets there, this latter topology is enough for the purpose of the
present paper, whenever considering sequences with continuous limit points.

If $f \in D$, we call \emph{local time process} of $f$ a Borel function
$(L(a,t), a, t \geq0)$ which satisfies
%
%
\begin{equation}
\label{eqnlocaltimedfn} \int_0^t \phi
\bigl(f(s)\bigr) \,ds = \int_0^\infty L(a,t) \phi(a)
\,da
\end{equation}
for any $t \geq0$ and any continuous function $\phi$ with compact
support included in $[0,\infty)$. The local time process of an
arbitrary function $f \in D$ may not exist, but when it does it is
unique (up to an almost everywhere modification). In the sequel we will
only consider local time processes associated to spectrally positive L\'
evy processes which either have infinite variation or negative drift.
These local time processes are known to exist; see, for instance,
Bertoin~\cite{Bertoin960}.

If $f \in D$ we define $\underline f$, the function $f$ reflected above
its past infimum, via
\[
\underline f(t) = f(t) - \min\Bigl(0, \inf_{0 \leq s \leq t} f(s)
\Bigr),\qquad
t \geq1.
\]

It is well known that this transformation induces a continuous map,
that is, if $f_n, f \in D$ are such that $f_n \to f$, then $\underline
{f}_n \to\underline f$.\vspace*{1.5pt}

For $f \in D$, let $\Delta f(t) = f(t) - f(t-)$ for $t > 0$, and $T_f =
\inf\{ t > 0\dvtx f(t) = 0 \}$ be the first time after time $0$ at
which~$f$ visits~$0$, with $T_f = \infty$ if $f$ never visits $0$ in
$(0,\infty)$. We see it as a map $T\dvtx D \to[0,\infty]$, and we
sometimes write $T(f)$ for~$T_f$. In general this map is not
continuous, but we have the following result.

%
\begin{lemma} \label{lemmaliminf}
If $f_n, f \in D$, $f_n \to f$ and $f$ is continuous, then $T_f
\leq\break
\liminf_n T_{f_n}$.
\end{lemma}

\begin{pf}
Let $\tau= \liminf_n T_{f_n}$ and $(u(n))$ such that $T_{f_{u(n)}}
\to\tau$. Since $f_n \to f$ and $f$ is continuous, we obtain
$f_{u(n)}(T_{f_{u(n)}}) \to f(\tau)$, hence $f(\tau) = 0$ which proves
the result.
\end{pf}

From now on $\Rightarrow$ denotes weak convergence. When considering
random vectors, we consider convergence in the product topology. The
previous lemma has the following consequence.

%
\begin{corollary} \label{corjoint}
If $X_n, X$ are stochastic processes such that $X_n \Rightarrow X$,
$X$ is continuous and $T_{X_n} \Rightarrow T_X$, then $(X_n, T_{X_n})
\Rightarrow(X, T_{X})$.
\end{corollary}

\begin{pf}
The sequence $(X_n, T_{X_n})$ is tight: let $(X', T')$ be any
accumulation point, so that $X'$ is equal in distribution to $X$ and
$T'$ to $T_X$. We show that $T' = T_{X'}$, which will show that $(X',
T')$ is equal in distribution to $(X, T_X)$ and will prove the result.
Assume without loss of generality that $(X_n, T_{X_n}) \Rightarrow(X',
T')$: the continuous mapping theorem and Lemma~\ref{lemmaliminf} imply
that $T_{X'} \leq T'$. But since they are equal in distribution, they
must be equal almost surely, hence the result.
\end{pf}

\subsection*{Stopping and shift operators}

For $t \geq0$ let $\sigma_t$ and $\theta_t$ be the stopping and shift
operators, respectively: for $f \in D$ and $t \geq0$, $\sigma_t f = f(
\cdot \wedge t)$ and $\theta_t f = f( \cdot + t)$. Note also
for simplicity $\sigma= \sigma_T$ and $\theta= \theta_T$, that is,
$\sigma f = \sigma_{T(f)} f$ and $\theta f = \theta_{T(f)} f$.
Formally, $\theta$ is only well defined if $T(f)$ is finite, and in the
rest of the paper we will only apply the map $\theta$ to such functions.

%
\begin{lemma} \label{lemmacontinuity-shift+stopping}
If $f_n, f \in D$ and $t_n, t \geq0$ are such that $f_n \to f$, $f$
is continuous and $t_n \to t$, then $\theta_{t_n} f_n \to\theta_{t} f$
and $\sigma_{t_n} f_n \to\sigma_{t} f$.
\end{lemma}

\begin{pf}
Let $w$ be the modulus of continuity of $f$, defined for $m,
\varepsilon> 0$ by $w_m(\varepsilon) = \sup\{ |f(t) - f(s)| \dvtx0
\leq
s, t \leq m \mbox{ and } |t-s| \leq\varepsilon\}$. Since $f$ is
continuous we have $w_m(\varepsilon) \to0$ as $\varepsilon\to0$, for
any $m \geq0$. Let $\overline t = \sup_{n \geq1} t_n$: for any $0
\leq s \leq m$, we have
\begin{eqnarray*}
\bigl\llvert\theta_{t_n} f_n(s) - \theta_{t}
f(s) \bigr\rrvert& \leq&\bigl\llvert f_n(s+t_n) -
f(s+t_n) \bigr\rrvert+ \bigl\llvert f(s+t_n) - f(s+t)
\bigr\rrvert
\\
& \leq&\norm{f_n - f}_{m + \overline t} + w_{m + \overline t}
\bigl(|t_n-t|\bigr)
\end{eqnarray*}
and similarly, $|\sigma_{t_n} f_n(s) - \sigma_{t} f(s)| \leq\norm{f_n
- f}_{m} + w_{m} (|t_n-t|)$. These upper bounds are uniform in $s \leq
m$, and since $f_n \to f$ and $t_n \to t$, letting $n \to+\infty$
gives the result.
\end{pf}

In the sequel we say that a sequence $(X_n)$ is C-tight if it is tight
and any accumulation point is almost surely continuous. We will use
several times that if $(X_n)$ and $(Y_n)$ are two C-tight sequences
defined on the same probability space, then the sequence $(X_n + Y_n)$
is also C-tight; see, for instance, Corollary VI.$3.33$ in Jacod and
Shiryaev~\cite{Jacod030}.

%
\begin{corollary} \label{cortightness-stopping+shift}
If $(X_n)$ is a C-tight sequence of processes and $(\kappa_n)$ is a
tight sequence of positive random variables, then $(\sigma_{\kappa_n}
X_n)$ and $(\theta_{\kappa_n} X_n)$ are C-tight.
\end{corollary}

\begin{pf}
Let $(u(n))$ be a subsequence, we must find $(v(n))$ a subsequence of
$(u(n))$ such that $(\sigma_{\kappa_{v(n)}} X_{v(n)})$ and $(\theta
_{\kappa_{v(n)}} X_{v(n)})$ converge weakly to a continuous process.
The sequence $(X_n, \kappa_n)$ being tight, there exists $(v(n))$ a
subsequence of $(u(n))$ such that $(X_{v(n)}, \kappa_{v(n)})$ converges
weakly to some $(X, \kappa)$, with $X$ a continuous process. Thus
$\sigma_{\kappa_{v(n)}} X_n \Rightarrow\sigma_{\kappa} X$ and
$\theta
_{\kappa_{v(n)}} X_n \Rightarrow\theta_{\kappa} X$ by Lemma~\ref
{lemmacontinuity-shift+stopping} together with the continuous mapping
theorem, hence the result.
\end{pf}

\subsection*{Excursions}

A function $e \in D_+$ will be called an excursion if $e(t) = 0$ for
some $t > 0$ implies $e(u) = 0$ for all $u \geq t$. Observe that
excursions are allowed to start at~$0$. Write $\Ecal$ for the set of
excursions with finite length $T_e$. Let also $\Ecal' \subset\Ecal$ be
the subset of excursions $e \in\Ecal$ such that $\int(1/e)$ is
finite, where from now on we write $\int_a^b f = \int_a^b f(t) \,dt$ and
$\int f = \int_0^{T_f} f$, for $f \in D$. When dealing with excursions
we will use the canonical notation for stochastic processes and write
$\epsilon$ for the canonical map.

For $f \in D$ and $\varepsilon> 0$, let $e_\varepsilon(f)$ be the
first excursion $e$ of $f$ away from $0$ that satisfies $T_e >
\varepsilon$, and let $g_\varepsilon(f) < d_\varepsilon(f)$ be its left
and right endpoints. Note that there need not be such an excursion, but
in the rest of the paper we will only apply the maps $e_\varepsilon$ to
functions $f$ such that $e_\varepsilon(f)$ is well defined for every
$\varepsilon> 0$. Also, note that by definition, we have
$e_\varepsilon= \sigma\circ\theta_{g_\varepsilon}$ in the sense that
for any $f \in D$,
\[
e_\varepsilon(f) = (\sigma\circ\theta_{g_\varepsilon(f)} ) (f) = (
\sigma_{d_\varepsilon(f) - g_\varepsilon(f)} \circ\theta
_{g_\varepsilon(f)} ) (f).
\]

\subsection*{Lamperti transformation}

We will call \emph{Lamperti transformation} the map $\Lcal\dvtx
\Ecal\to
\Ecal$ that to an excursion $f \in\Ecal$ associates the excursion $h
\in\Ecal$ defined by $h(\int_0^t f) = f(t)$ for all $t \geq0$. More
specifically, if $\kappa$ is the inverse of the strictly increasing,
continuous function $t\mapsto\int_0^t f$ on $[0, \int f]$, then
$\Lcal
(f)= f\circ\kappa$ on $[0, \int f]$ and 0 elsewhere. In particular,
$(T\circ\Lcal)(f) = \int f$.

The inverse Lamperti transformation $\Lcal^{-1}$ also plays a crucial
role. By definition, $\Lcal^{-1}(f)$ is the solution in $\Ecal$, when
it exists and is unique, to the equation $h(t) = f(\int_0^t h)$, $t
\geq0$, where $h$ is the unknown function. Existence and uniqueness to
such equations are studied in Chapter~$6$ of Ethier and Kurtz~\cite
{Ethier860}. Because we consider excursions which may start at $0$, we
cannot directly invoke Theorem~$1.1$ there, but an inspection of the
proof reveals that it can be adapted to show that $\Lcal^{-1}(f)$ is
well defined for $f \in\Ecal'$. In this case, we have $\Lcal
^{-1}(f)=f\circ\pi$ on $[0, \int(1/f)]$, and $0$ otherwise, with
$\pi$
the inverse of the strictly increasing, continuous function $t\mapsto
\int_0^{t} (1/f)$ on $[0, \int(1/f)]$. In particular, $(T\circ\Lcal
^{-1})(f) = \int(1/f)$.

We will need the following results on $\Lcal$ and $\Lcal^{-1}$, which
are closely related to results by Helland~\cite{Helland780} or Ethier
and Kurtz~\cite{Ethier860}, Chapter~$6$. There are nonetheless
significant differences and for completeness, we provide the proof of
the following lemma in the \hyperref[app]{Appendix}.

%
\begin{lemma} \label{lemmacontinuity+tightness-lamperti}
Let $X_n, X$ be random elements of $\Ecal$ such that the sequence
$(X_n)$ is C-tight and the sequence $(T_{X_n})$ is tight.

If $X_n \Rightarrow X$ then $\Lcal(X_n) \Rightarrow\Lcal(X)$.

If $\P(X_n \in\Ecal') = 1$, then the sequence $(\Lcal^{-1}(X_n))$
is C-tight.

If $X_n \Rightarrow X$ and $\P(X_n \in\Ecal') = \P(X \in\Ecal') =
1$, then $\Lcal^{-1}(X_n) \Rightarrow\Lcal^{-1}(X)$.
\end{lemma}

\section{CMJ branching processes, PS queues and L\'evy processes}
\label{secCMJ-PS-Levy}

Recall from the \hyperref[sec1]{Introduction} that a
Crump--Mode--Jagers (CMJ) process
is a stochastic process with nonnegative integer values counting the
size of a population where individuals give birth to independent copies
of themselves, and that processor-sharing (PS) is the service
discipline where the server splits its service capacity equally among
all users present in the queue at any time.\looseness=1

In the sequel, \emph{we will only consider homogeneous and binary CMJ
processes}, where individuals give birth to a single offspring at times
of a Poisson process independent of their life length. With the
notation of the \hyperref[sec1]{Introduction}, $\xi$ is a Poisson
process independent of
$V$. Similarly, \emph{we will only consider} $M/G/1$-PS \emph{queues},
that is, PS queues with Poisson arrivals and i.i.d. service
requirements. Henceforth, CMJ will stand for homogeneous and binary
CMJ, and PS for $M/G/1$-PS.

In particular, thanks to the memoryless property of the exponential
random variable, both a CMJ process and a PS queue can be described by
a Markov process living in the state space $\Scal= \bigcup_{n \geq0}
(0,\infty)^n$, the set of finite sequences of positive real numbers.
For a CMJ process the Markovian descriptor keeps track of the residual
life lengths of the individuals alive; for a PS queue, it keeps track
of the residual service requirements of the customers present in the
queue. Although for the PS queue we will sometimes refer to this Markov
process, we will actually only consider marginals of it, namely the
workload process (corresponding to its total mass) and the queue length
process (corresponding to the cardinal of its support). Thus although
studying non-Markovian processes, we avoid the framework of measure
valued processes.

\subsection*{Simple facts about the PS queue}

Let $q$ be the queue length process of a PS queue with unit service
capacity, arrival rate $\lambda$ and service distribution~$S$. The
workload process is the process keeping track of the total amount of
work in the system, defined as the sum over all the customers of their
residual service requirements. Since we assume Poisson arrivals, the
workload process is a compensated compound Poisson process with
drift~$-1$ and L\'evy measure $\lambda\P(S\in\cdot)$, reflected above
its past infimum.

Set $\rho:=\lambda\E(S)$ the load, and assume $\rho<1$ (subcritical
case). Let $S^*$ be the random variable with density $\P(S \geq\cdot
) / \E(S)$ with respect to Lebesgue measure. It is sometimes called the
forward recurrence time of $S$ and has mean $\E(S^*) = \E(S^2) / (2\E
(S))$. The assumption $\rho< 1$ is equivalent to assuming that the
Markov process describing the PS queue has a unique invariant
distribution $\nu^*$ on $\Scal$. In that case, the invariant
distribution is characterized by a geometrically distributed number of
customers with parameter $\rho$ and i.i.d. residual service times with
common distribution~$S^*$; see, for example, Robert~\cite{Robert030},
Proposition~$7.13$.

\subsection*{Connection between PS queues, CMJ processes and L\'evy processes}

In the following statement, $q$ is the above PS queue, and $\P^\chi$ is
its law started at $\chi\in\Scal$. The following result is known
since at least Kitayev and Yashkov~\cite{Kitayev790}; see also
Chapter~$7.3$ in Robert~\cite{Robert030}.

%
\begin{theorem}[(Connection between PS queues and CMJ processes)]\label{thmtime-change}
Let $\chi= (\chi_i, 1 \leq i \leq k) \in\Scal$. The process $\Lcal
^{-1}(q)$ under $\P^\chi$ is a CMJ process starting with~$k$ ancestors,
with birth rate $\lambda$ and life length distribution $S$, except for
the ancestors who have deterministic life lengths given by~$\chi$.
\end{theorem}

Thus we can see $\sigma q$ as the time change of a CMJ process, since
$\sigma q = \Lcal(z)$ with $z = \Lcal^{-1}(q)$ which by the above is a
CMJ process. Further, since $0$ is a regeneration point of $q$, every
excursion of $q$ away from $0$ can be seen as the time change of a CMJ
process started with one individual with life length distributed as $S$.

The following result can be found in Lambert~\cite{Lambert100}. The
jumping contour process of a homogeneous, binary CMJ tree starting from
one progenitor is the key object underlying this result. It is defined
in Lambert~\cite{Lambert100}, to which the reader is referred for more details.

In the following statement, $x$ denotes a spectrally positive L\'evy
process starting from $\delta>0$, with drift $-1$, L\'evy measure
$\lambda\P(S\in\cdot)$ and local time process $(\ell(a,t), a, t
\geq
0)$, as defined in \eqref{eqnlocaltimedfn}. Note that in this case,
$\ell(a,t)$ is also the number of times when $x$ has taken the value
$a$ before time $t$.

%
\begin{theorem}[(Connection between CMJ and L\'evy processes)]\label
{thmLambert}
The process $(\ell(a,T_x), a \geq0)$ is a CMJ process with birth rate
$\lambda$ and life length distribution~$S$, started with one progenitor
with life length $\delta$.
\end{theorem}

\subsection*{Scaling near the critical point}

For each integer $n\ge1$, consider some $\lambda_n>0$ and a positive
random variable $S_n$, with forward recurrence time $S_n^*$. Let $q_n$
denote the queue length process of a PS queue with arrival
rate~$\lambda
_n$ and service distribution~$S_n$. Let $z_n = \Lcal^{-1}(q_n)$, which
according to Theorem~\ref{thmtime-change} is a CMJ process with birth
rate~$\lambda_n$ and life length distribution $S_n$. Last, let $x_n$ be
a compensated compound Poisson process with drift $-1$ and L\'evy
measure $\lambda_n\P(S_n\in\cdot)$. Let $(\ell_n(a,t), a, t \geq0)$
be its local time process, so that by Theorem~\ref{thmLambert}, $z_n$
is equal in distribution to $(\ell_n(a, T_{x_n}), a \geq0)$. Also,
$x_n$ is equal in distribution to the workload process corresponding to
$q_n$ (with suitable initial conditions); in particular, the zero sets
of $q_n$ and $x_n$ have the same distribution.

We restrict our attention to the subcritical case; namely, we assume
that for each $n \geq1$ the load $\rho_n:= \lambda_n \E(S_n)$
satisfies $\rho_n < 1$. This assumption means that all hitting times
of~$0$ by $q_n$, $z_n$ and $x_n$ with deterministic initial states
[$x_n$ starting in $(0,\infty)$] are almost surely finite and have
finite expectations. We consider the following scaling near the
critical point: in the sequel we assume that there exist finite and
strictly positive real numbers $\lambda, \beta$ and~$\alpha$ such that
%
%
\begin{equation}
\label{eqHT} \qquad\lim_{n \to+\infty} \lambda_n = \lambda,\qquad \lim
_{n \to+\infty} \frac
{\lambda_n}{2}\E\bigl(S_n^2
\bigr) = \beta\quad\mbox{and}\quad \lim_{n \to+\infty} n (1-\rho_n)
= \alpha.
\end{equation}

These three assumptions imply that $\E(S_n^*) \to\beta$. We are
interested in the rescaled processes $Q_n$, $Z_n$ and $X_n$ defined by
\begin{eqnarray}
Q_n(t) = \frac{q_n(n^2 t)}{n},\qquad Z_n (t) = \frac{z_n ( n t )}{n}\quad
\mbox{and}\quad X_n(t) = \frac{x_n(n^2 t)}{n},\nonumber\\
  \eqntext{n \geq1, t \geq0.}
\end{eqnarray}

Let also
\[
L_n(a,t) = \frac{\ell_n(na, n^2t)}{n},\qquad  n \geq1, a, t \geq0.
\]

Then $L_n$ is the local time process of $X_n$, that is, it satisfies
\[
\int_0^t \phi\bigl(X_n(s)\bigr)
\,ds = \int_0^\infty\phi(a) L_n(a,t)
\,da.
\]

By the L\'evy--Khintchine formula, $X_n$ is a L\'evy process with
Laplace exponent $\Psi_n(u) = n u - n^2 \lambda_n \E(1 - e^{-u S_n /
n})$. In view of~\eqref{eqHT}, standard arguments show that $\Psi_n(u)
\to\alpha u + \beta u^2$ for any $u \geq0$. As a consequence, see for
instance Kallenberg~\cite{Kallenberg020}, $(X_n)$ converges in
distribution to a drifted Brownian motion with drift $-\alpha$ and
Gaussian coefficient $2\beta$, which we write in the sequel $X$ and
whose local time process is denoted $(L(a,t), a, t \geq0)$.

\subsection*{Notation for the initial condition}

For any $\chi=(\chi_1,\ldots,\chi_k)\in\Scal$, when $q_n$ and
$z_n$ are
started with $k \geq1$ customers/individuals with residual service
times/life lengths given by~$\chi$, the probability measure is denoted
$\P_n^\chi$. The law of the PS queue started empty will be
denoted~$\P
_n^\varnothing$.
We will use the following notation for random initial conditions.

When~$q_n$ and~$z_n$ are started with \emph{one} individual with
residual life length distributed as~$S_n$, the law will simply be
denoted~$\P_n$. When there are initially~$\zeta_n$ individuals with
i.i.d. residual life lengths distributed as~$S_n^*$, we will use the
symbol~$\P_n^{\zeta_n *}$. When the initial condition is~$\nu_n^*$ (a
geometric number with parameter~$\rho_n$ of individuals with i.i.d.
life lengths distributed as~$S_n^*$), we will merely use the symbol~$\P
_n^{*}$. Note that $Q_n$ under $\P_n^*$ is a stationary process.

The probability measure for the L\'evy processes is denoted $\PX_n^a$
when $X_n$ itself is started at $a \in\R$. When $X_n$ starts at a
random initial value distributed as $S_n/n$ we write $\PX_n$. Finally,
$\PX^a$ denotes the law of $X$ started at~$a$.

The scalings in time and space leading to $Q_n$, $Z_n$ and $X_n$ have
been chosen in order to preserve the defining relationships between
$q_n$, $z_n$ and $x_n$.

%
\begin{lemma} \label{lemmarescaling}
We have $Z_n = \Lcal^{-1}(Q_n)$ and $\sigma Q_n = \Lcal(Z_n)$, in
particular $T_{Q_n} = \int Z_n$. Moreover, for any $\delta> 0$, $Z_n$
under $\P_n^{\chi}$ with $\chi= (n\delta) \in\Scal$ is equal in
distribution to $(L_n(a, T_{X_n}), a \geq0)$ under $\PX^{\delta}_n$.
\end{lemma}

\subsection*{Excursion measures}

In the sequel, three distinct excursion measures will be considered.
First, $\Ncal$ is the excursion measure of $\underline X$ away from 0,
where in this case, the excursion measure is normalized so that the
local time of $\underline{X}$ at 0 at time $t$ is taken equal to
$-\min
(0, \inf_{0\le s\le t} X(t))$. This normalization will always be
considered for processes reflected above their past infimum.

Second, $\Mcal$ the push-forward of $\Ncal$ by the map $L( \cdot,
T) = (L(a, T), a \geq0)$, where from now on we will also denote by
$(L(a,t), a, t \geq0)$ the local time process of the canonical
excursion $\epsilon$. In other words, for any measurable function
$f\dvtx\Ecal\to[0,\infty)$, we have
\[
\Mcal(f ) = \Ncal\bigl(f \circ L( \cdot, T) \bigr).
\]

Third, we define $\Ncal'$ as the push-forward of $\Mcal$ by $\Lcal$,
that is, for any measurable function $f\dvtx\Ecal\to[0,\infty)$, we have
\[
\Ncal' (f ) = \Mcal(f \circ\Lcal).
\]

Then the measures obtained by taking the push-forward of $\Ncal$ and
$\Ncal'$ by $T$ coincide, that is, for any Borel set $A \subset
[0,\infty)$ we have
%
%
\begin{equation}
\label{eqequality} \Ncal(T \in A) = \Ncal'(T \in A).
\end{equation}

Indeed, we have by definition of $\Ncal'$, $\Mcal$ and $L$,
\[
\Ncal' (T \in A ) = \Mcal\biggl(\int\epsilon\in A \biggr) = \Ncal
\biggl(\int_0^\infty L(a,T) \,da \in A \biggr) =
\Ncal(T \in A)
\]
since $\int_0^\infty L(a, T) \,da = T$. As a side remark, note that it
could be proved that $\Mcal(1 \wedge T) = +\infty$, and so there is no
regenerative process admitting $\Mcal$ as its excursion measure.

\section{Heavy traffic of PS via excursion theory} \label{secmain}

The goal of this section is to prove forthcoming Theorem~\ref{thmconv-PS}.
Roughly speaking, it states that the sequence $(Q_n)$ of PS
queue length processes started empty converge weakly to the
regenerative process with excursion measure away from $0$ equal to
$\Ncal'$. Recall that $\Ncal'$ is the push-forward of the excursion
measure $\Ncal$ of $\underline{X}$ by the successive application of
$L(\cdot,T)$ (local time process at the first hitting time of 0) and
$\Lcal$ (Lamperti transformation).

The push-forward $\Mcal$ of $\Ncal$ by the mere application of
$L(\cdot,T)$ is not an excursion measure [in the above mentioned sense that
$\Mcal(1 \wedge T) = + \infty$], but we expect nonetheless that the
distributions of the CMJ processes $Z_n$ will converge in some sense to
$\Mcal$. This intuition is made precise in Theorem~\ref
{thmconv-CMJ-excursion}, where $Z_n$ starts with one initial
individual and is suitable conditioned, and Theorem~\ref
{thmconv-CMJ-star}, where $Z_n$ starts from a large initial condition.

We also specify that $\Ncal'$ will be identified in Theorem~\ref
{thmidentification} as the excursion measure away from 0 of $\beta
^{-1}\underline{X}$, which is the reflected Brownian motion with drift
$-\alpha/\beta$ and Gaussian coefficient $2/\beta$. This will ensure
that the sequence $(Q_n)$ actually converges weakly to this reflected
process; see also Theorem~\ref{thmconv-PS-2} for general initial condition.

%
\begin{theorem} \label{thmconv-PS}
Let $Q_\infty$ be the process obtained by applying It\^o's
construction to the excursion measure $\Ncal'$. Then the sequence
$(Q_n)$ under $\P_n^\varnothing$ converges weakly to $Q_\infty$.
\end{theorem}

To avoid any ambiguity, let us explain what we mean by It\^o's
construction; see Blumenthal~\cite{Blumenthal920}, for instance. Let
$\partial$ be some cemetery point and $e = (e_t, t \geq0)$ be an
$\Ecal\cup\{\partial\}$-valued Poisson point process with intensity
measure~$\Ncal'$. Define
\[
\widetilde L(t) = \sum_{0 \leq s \leq t} T(e_s)
\]
with the convention $T(\partial) = 0$. Since $\Ncal'(1 \wedge T) <
+\infty$ by~\eqref{eqequality}, $\widetilde L$ is a subordinator with
L\'evy measure $\Ncal'(T \in\cdot )$. Let $\widetilde L^{-1}$ be
the right-continuous inverse of $\widetilde L$; then the process
$Q_\infty$ is defined via the following formula:
\[
Q_\infty(t) = e_{\widetilde L^{-1}(t-)} \bigl(t - \widetilde L\bigl
(\widetilde
L^{-1}(t)-\bigr) \bigr) \mathbh{1}_{\Delta\widetilde
L(\widetilde L^{-1}(t))\not=0},\qquad t \geq0.
\]

We first prove in Section~\ref{sublevy} preliminary results on L\'evy
processes. Section~\ref{subtightness} is devoted to tightness,
Section~\ref{subCMJ} proves a result of independent interest on CMJ
processes and Section~\ref{subproof} provides the proof of
Theorem~\ref
{thmconv-PS}.

\subsection{Preliminary results on L\'evy processes} \label{sublevy}

We will need the following results on L\'evy processes.

%
\begin{lemma} \label{lemmalevy-T}
For any $a > 0$, the sequence $(X_n, T_{X_n})$ under $\PX_n^a$
converges weakly to $(X, T_X)$ under~$\PX^a$.
\end{lemma}

\begin{pf}
Since $\PX^a(\forall\varepsilon> 0, \inf_{[0,\varepsilon]} \theta
X < 0) = 1$, the result follows directly from Proposition~VI.$2.11$ in
Jacod and Shiryaev~\cite{Jacod030}.
\end{pf}

%
\begin{lemma} \label{lemmalevy-endpoints}
For any $\varepsilon> 0$, the sequence $(g_\varepsilon(\underline
X_n), d_\varepsilon(\underline X_n))$ under $\PX_n^0$ converges weakly
to $(g_\varepsilon(\underline X), d_\varepsilon(\underline X))$ under
$\PX^0$.
\end{lemma}

\begin{pf}
Remember that $\Psi_n$ is the Laplace exponent of $X_n$. Since $X_n$
drifts to $-\infty$, $\Psi_n$ is continuous and strictly increasing and
we denote $\Phi_n$ its inverse. Let similarly $\Psi$ be the Laplace
exponent of $X$ and $\Phi$ its inverse. Since $X_n \Rightarrow X$ it is
not hard to show that $\Phi_n(u) \to\Phi(u)$ for every $u \geq0$.
Moreover, for $t \geq0$ let
\[
\gamma_n(t) = \inf\bigl\{ s \geq0\dvtx X_n(s) = -t
\bigr\} \quad\mbox{and}\quad \gamma(t) = \inf\bigl\{ s \geq0\dvtx X(s) =
-t \bigr\}.
\]
Then it is well known (see, e.g., Bertoin~\cite{Bertoin960},
Theorem~VII.$1$) that $\gamma_n$ and $\gamma$ are
subordinators with Laplace exponent $\Phi_n$ and $\Phi$, respectively.
Since $\Phi_n(u) \to\Phi(u)$ for every $u \geq0$, standard arguments
imply that $\gamma_n \Rightarrow\gamma$. Moreover, since $\gamma_n$
and $\gamma$ are the right-continuous inverses of the local time
processes of $\underline X_n$ and $\underline X$ at $0$, we have the
identities $g_\varepsilon(\underline X_n) = \gamma_n(t^1_\varepsilon
(\gamma_n)-)$ and $d_\varepsilon(\underline X_n) = \gamma
_n(t^1_\varepsilon(\gamma_n))$ and similarly without the subscript $n$,
where in the rest of the proof we define $t^1_\varepsilon(f) = \inf\{t
\geq0\dvtx|\Delta f(t)| > \varepsilon\}$ for any $f \in D$ and
$\varepsilon> 0$.

Proposition~$2.7$ in Jacod and Shiryaev~\cite{Jacod030} shows that if
$f_n, f \in D$ are such that $f_n \to f$, $t^1_\varepsilon(f) <
+\infty
$ and $\varepsilon\notin\{ |\Delta f(t)|\dvtx t \geq0 \}$, then
$f_n(t^1_\varepsilon(f_n)-) \to f(t^1_\varepsilon(f)-)$ as well as
$f_n(t^1_\varepsilon(f_n)) \to f(t^1_\varepsilon(f))$. The desired
result therefore follows from an application of the continuous mapping
theorem, together with the fact that $\PX^0(t^1_\varepsilon(\gamma) <
+\infty, \varepsilon\notin\{ |\Delta\gamma(t)|\dvtx t \geq0 \}) =
1$ for
every $\varepsilon> 0$.
\end{pf}

%
\begin{lemma}\label{lemmalevy-exc}
For any $\varepsilon> 0$, the sequence $(\sigma X_n, T_{X_n})$
considered under $\PX_n( \cdot | T_{X_n} > \varepsilon)$
converges weakly to $(\epsilon, T_\epsilon)$ under~$\Ncal( \cdot |
T > \varepsilon)$.
\end{lemma}

\begin{pf}
From now on and unless otherwise specified, we implicitly consider
$X_n$ under $\PX_n^0$ and $X$ under $\PX^0$. Since by definition the
process $\sigma X_n$ under $\PX_n( \cdot | T_{X_n} > \varepsilon
)$ is equal in distribution to $e_\varepsilon(\underline X_n)$ and
$\Ncal( \cdot | T > \varepsilon)$ is the law of $e_\varepsilon
(\underline X)$, the result is equivalent to showing that
$(e_\varepsilon, T \circ e_\varepsilon)(\underline X_n) \Rightarrow
(e_\varepsilon, T \circ e_\varepsilon)(\underline X)$.

For $f \in D$ let $\Jcal(f) = (\underline f, g_\varepsilon(\underline
f), d_\varepsilon(\underline f))$. In view of Lemma~\ref
{lemmacontinuity-shift+stopping} and the continuous mapping theorem,
to prove that $(e_\varepsilon, T \circ e_\varepsilon)(\underline X_n)
\Rightarrow(e_\varepsilon, T \circ e_\varepsilon)(\underline X)$ it is
enough to show that $\Jcal(X_n) \Rightarrow\Jcal(X)$. We have
$\underline X_n \Rightarrow\underline X$, while Lemma~\ref
{lemmalevy-endpoints} shows that $(g_\varepsilon, d_\varepsilon
)(\underline X_n) \Rightarrow(g_\varepsilon, d_\varepsilon
)(\underline
X)$. Hence the sequence $(\Jcal(X_n))$ is tight, and we only need to
identify accumulation points. Let $(X', g', d')$ be any accumulation
point, and assume without loss of generality using Skorohod's
representation theorem that $\Jcal(X_n) \to(X', g', d')$: then $X'$ is
equal in distribution to $\underline X$ and $(g', d')$ to
$(g_\varepsilon, d_\varepsilon)(\underline X)$, and we only have to
show that $(g', d') = (g_\varepsilon, d_\varepsilon)(X')$.

Since $\underline X_n$ and $X_n$ have the same generator in $(0,\infty
)$ and the functional $T$ is $\PX$-a.s. continuous (see Lemma~\ref
{lemmalevy-T}), it can be proved that $T(\theta_{t_n} \underline X_n)
\to T(\theta_t X')$ for any $t_n, t \geq0$ such that $t_n \to t$.

Let $t < g'$: since $g_\varepsilon(\underline X_n) \to g'$ we have $t
< g_\varepsilon(\underline X_n)$ for $n$ large enough, and for those
$n$s it holds by definition of $g_\varepsilon(\underline X_n)$ that
$T(\theta_t \underline X_n) \leq\varepsilon$. Since $T(\theta_{t}
\underline X_n) \to T(\theta_t X')$ we obtain that $T(\theta_t X')
\leq
\varepsilon$. Since $t < g'$ is arbitrary, this proves that
$g_\varepsilon(X') \geq g'$, and since they are equal in distribution,
they must be equal almost surely. Since $T(\theta_{g_\varepsilon
(\underline X_n)}\underline X_n) = d_\varepsilon(\underline X_n) -
g_\varepsilon(\underline X_n)$, letting $n \to+\infty$ shows that
$T(\theta_{g'} X') = d' - g'$, and so $d' = g_\varepsilon(X') +
T(\theta
_{g_\varepsilon(X')}X') = d_\varepsilon(X')$. The proof is complete.
\end{pf}

\subsection{Tightness} \label{subtightness}

Although tightness is usually a technical issue, it comes here from a
simple queueing argument. Theorem~\ref{thmdeparture-process} may look
naive to an experienced reader, but to the best of our knowledge its
implications in terms of tightness have never been used before; similar
arguments could, for instance, have been used in Limic~\cite
{Limic010}. Since processor-sharing is a symmetric service discipline,
the following result is a direct consequence of Theorems~$3.10$
and~$3.6$ in Kelly~\cite{Kelly790}.

%
\begin{theorem} \label{thmdeparture-process}
The departure process of the queue length process $q_n$ under $\P_n^*$
is a Poisson process with parameter~$\lambda_n$.
\end{theorem}

%
\begin{corollary} \label{cortightness-stationary-Q}
The sequence of processes $(Q_n)$ under $\P_n^*$ is C-tight.
\end{corollary}

\begin{pf}
Writing $a_n$ and $d_n$ for the arrival and departure processes,
respectively, we can write $Q_n(t) = Q_n(0) + A_n(t) - D_n(t)$ for $t
\geq0$ with $A_n(t) = (a_n(n^2 t) - n^2 \lambda_n t) / n$ and $D_n(t)
= (d_n(n^2 t) - n^2 \lambda_n t) / n$. By Theorem~\ref
{thmdeparture-process}, $a_n$ and $d_n$ under $\P_n^*$ are two Poisson
processes with intensity~$\lambda_n$. Since $\lambda_n \to\lambda$,
both $(A_n)$ and $(D_n)$ converge in distribution to a Brownian motion
and so are C-tight, and hence so is the difference $(A_n - D_n)$. Since
$(Q_n(0))$ under $\P_n^*$ converges to an exponential random variable,
this shows that $(Q_n)$ under $\P_n^*$ is C-tight.
\end{pf}

Corollary~\ref{cortightness-stationary-Q} encompasses all the
tightness results we need. To be more specific, by considering $Q_n$
under $\P_n^*$ and shifting it at time $T_{Q_n}$ we can get the
tightness of $(Q_n)$ under $\P_n^\varnothing$. Also, we can get the
tightness of CMJ processes by suitably selecting excursions of $Q_n$
and applying $\Lcal^{-1}$. The elementary operations that we need to
perform preserve C-tightness by Corollary~\ref
{cortightness-stopping+shift} and Lemma~\ref
{lemmacontinuity+tightness-lamperti}, and so we get the following result.

%
\begin{corollary} \label{cortightness}
The sequence $(Q_n)$ under $\P_n^\varnothing$ is C-tight, and for any
$\varepsilon> 0$, the sequence $(Z_n)$ under $\P_n ( \cdot |
\int Z_n > \varepsilon)$ is C-tight.
\end{corollary}

\begin{pf}
By regeneration of $Q_n$, $Q_n$ under $\P_n^\varnothing$ is equal in
distribution to $\theta Q_n$ under $\P_n^*$. Since $T_{Q_n}$ is equal
in distribution to $T_{X_n}$, Lemma~\ref{lemmalevy-T} shows that
$(T_{Q_n})$ under $\P_n^*$ is tight. Combining Corollaries~\ref
{cortightness-stationary-Q} and~\ref{cortightness-stopping+shift}
shows the C-tightness of $(Q_n)$ under $\P_n^\varnothing$.

For $Z_n$, since $Z_n = \Lcal^{-1}(Q_n)$ and $\int Z_n = T_{Q_n}$ by
Lemma~\ref{lemmarescaling}, one sees that $Z_n$ under $\P_n( \cdot
| \int Z_n > \varepsilon)$ is equal in distribution to $\Lcal
^{-1}(e_\varepsilon(Q_n))$ under $\P_n^\varnothing$. Since the zero set\vadjust{\goodbreak}
of $Q_n$ is equal in distribution to the zero set of $\underline X_n$,
$(g_\varepsilon, d_\varepsilon)(Q_n)$ under $\P_n^\varnothing$ is equal
in distribution to $(g_\varepsilon, d_\varepsilon)(\underline X_n)$
under $\PX_n^0$. Thus Lemma~\ref{lemmalevy-endpoints} implies that the
sequence $(g_\varepsilon, d_\varepsilon)(Q_n)$ under $\P
_n^\varnothing$
is tight. Since by definition $e_\varepsilon(Q_n) = (\sigma
_{d_\varepsilon(Q_n) - g_\varepsilon(Q_n)} \circ\theta
_{g_\varepsilon
(Q_n)}) (Q_n)$, combining the results of Corollary~\ref
{cortightness-stopping+shift} and Lemma~\ref
{lemmacontinuity+tightness-lamperti} and using also that $(Q_n)$ under
$\P_n^\varnothing$ is C-tight, we obtain the C-tightness of $(Z_n)$ under
$\P_n( \cdot | \int Z_n > \varepsilon)$.
\end{pf}

\subsection{Weak convergence of CMJ processes} \label{subCMJ}

The following result is of independent interest in the area of
branching processes. Using similar techniques and ideas, the
conditionings $\{ \int Z_n > \varepsilon\}$ and $\{ \int\epsilon>
\varepsilon\}$ in the next statement could be replaced by $\{ T_{Z_n}
> \varepsilon\}$ and $\{ T_{\epsilon} > \varepsilon\}$, respectively.
Theorem~\ref{thmconv-CMJ-star} in the following section gives another
result with a large initial condition.

%
\begin{theorem} \label{thmconv-CMJ-excursion}
For any $\varepsilon> 0$, the sequence $(Z_n)$ under $\P_n ( \cdot
| \int Z_n > \varepsilon)$ converges weakly to $\Mcal( \cdot|
\int\epsilon> \varepsilon)$.
\end{theorem}

\begin{pf}
In the remainder of the proof we implicitly consider $Z_n$ under $\P_n
( \cdot | \int Z_n > \varepsilon)$ and $X_n$ under $\PX_n(
\cdot | T_{X_n} > \varepsilon)$ and we denote by $L_n^0$ the
process $(L_n(a, T_{X_n}), a \geq0)$. Lemma~\ref{lemmarescaling}
shows that $Z_n$ is equal in distribution to $L_n^0$, so we only have
to show that $L_n^0 \Rightarrow\Mcal( \cdot | \int\epsilon>
\varepsilon)$.

Lemma~\ref{cortightness} shows that the sequence $(L_n^0)$ is
C-tight, so we only have to identify accumulation points. Let $Z$ be
any accumulation point and assume without loss of generality that
$L_n^0 \Rightarrow Z$. Lemma~\ref{lemmalevy-exc} shows that $(\sigma
X_n, T_{X_n})$ converges weakly to $(\epsilon, T_\epsilon)$ under
$\Ncal
( \cdot | T > \varepsilon)$. Then the sequence $(\sigma X_n,
T_{X_n}, L_n^0)$ is tight. Let $(e, \tau, Z')$ be any accumulation
point, so that $(e, \tau)$ is equal in distribution to $(\epsilon,
T_\epsilon)$ under $\Ncal( \cdot | T > \varepsilon)$ and $Z'$
to $Z$. Assume without loss of generality by Skorohod's representation
theorem that $(\sigma X_n, T_{X_n}, L_n^0) \to(e, \tau, Z')$. By
definition we have
\[
\int_0^{T_{X_n}} \phi\bigl(\sigma X_n(t)
\bigr) \,dt = \int_0^\infty\phi(a)
L_n^0(a) \,da
\]
for all continuous functions $\phi$ with a compact support. Thus
passing to the limit, the dominated convergence theorem (or uniform
convergence arguments) shows that
\[
\int_0^{\tau} \phi\bigl(e(t)\bigr) \,dt = \int
_0^\infty\phi(a) Z'(a) \,da,
\]
which shows that $Z'$ is the local time process of $e$ up to time $\tau
$. The result is proved.
\end{pf}

\subsection{\texorpdfstring{Proof of Theorem~\protect\ref{thmconv-PS}}{Proof of Theorem 4.1}} \label{subproof}

We begin with a preliminary result.

%
\begin{lemma}\label{lemmaexc}
For any $\varepsilon> 0$, the sequence $(\sigma Q_n, T_{Q_n})$
considered under $\P_n( \cdot | T_{Q_n} > \varepsilon)$
converges weakly to $(\epsilon, T_\epsilon)$ under $\Ncal'( \cdot
| T > \varepsilon)$.
\end{lemma}

\begin{pf}
Until the end of this step, we consider implicitly the process $Q_n$,
and hence $Z_n$, under $\P_n( \cdot | T_{Q_n} > \varepsilon) =
\P_n( \cdot | \int Z_n > \varepsilon)$. By Theorem~\ref
{thmconv-CMJ-excursion}, we know that $Z_n \Rightarrow\Mcal( \cdot
| \int\epsilon> \varepsilon)$. Moreover, Lemma~\ref
{lemmarescaling} implies that $T_{Z_n}$ is equal in distribution to
$\norm{X_n}_{T_{X_n}}$ under $\PX_n( \cdot | T_{X_n} >
\varepsilon)$ which converges, in view of Lemma~\ref{lemmalevy-exc}
and using the continuous mapping theorem, to $\norm{\epsilon}_\infty$
under $\Ncal( \cdot | T > \varepsilon)$. In particular, the
sequence $(T_{Z_n})$ is tight, so Lemma~\ref
{lemmacontinuity+tightness-lamperti} implies that the sequence $(\Lcal
(Z_n))$ converges weakly to the push-forward of $\Mcal( \cdot |
\int\epsilon> \varepsilon)$ by $\Lcal$. Since $\Lcal(Z_n) = \sigma
Q_n$ by Lemma~\ref{lemmarescaling} and the push-forward of $\Mcal(
\cdot | \int\epsilon> \varepsilon)$ by $\Lcal$ is by definition
equal to $\Ncal'( \cdot | T > \varepsilon)$, and we obtain the
convergence of the sequence $(\sigma Q_n)$ toward $\Ncal'( \cdot
| T > \varepsilon)$.

On the other hand, since the workload associated to $Q_n$ is equal in
distribution to $\underline X_n$, we obtain that $T_{Q_n}$ is equal in
distribution to $T_{X_n}$ under $\PX_n( \cdot | T_{X_n} >
\varepsilon)$, hence $(T_{Q_n})$ converges weakly to $T$ under $\Ncal
( \cdot | T > \varepsilon)$ in view of Lemma~\ref
{lemmalevy-exc}. Since $T$ under $\Ncal( \cdot | T >
\varepsilon)$ is equal in distribution to $T$ under $\Ncal'( \cdot
| T > \varepsilon)$ by~\eqref{eqequality}, and we obtain the
convergence of $(T_{Q_n})$ toward $T$ under $\Ncal'( \cdot | T
> \varepsilon)$. To conclude that the joint convergence holds we invoke
Corollary~\ref{corjoint}.
\end{pf}

We now prove Theorem~\ref{thmconv-PS}. Since the sequence $(Q_n)$
under $\P_n^\varnothing$ is C-tight by Corollary~\ref{cortightness}, we
only have to identify accumulation points. So let $Q$ be any
accumulation point, and assume without loss of generality that $Q_n
\Rightarrow Q$: we must prove that $Q$ is equal in distribution to
$Q_\infty$. In the rest of this section, for $\varepsilon> 0$ let
$A_\varepsilon\dvtx D \to\Ecal\times[0,\infty) \times[0,\infty
)$ be the
map given by $A_\varepsilon= (e_\varepsilon, g_\varepsilon,
d_\varepsilon)$, and let $\Phi_\varepsilon\dvtx D \to D$ the map that
truncates excursions with length smaller than $\varepsilon$; that is,
for $f \in D$ and $t \geq0$ we put $\Phi_\varepsilon(f)(t) = f(t)$ if
$f(t) \neq0$ and the excursion $e$ of $f$ straddling $t$ satisfies
$T_e > \varepsilon$; otherwise we put $\Phi_\varepsilon(f)(t) = 0$. We
prove that $Q$ is equal in distribution to $Q_\infty$ in two
steps.\looseness=-1

\textit{First step}. Let $\varepsilon> 0$: we first prove
that $(Q_n, A_\varepsilon(Q_n)) \Rightarrow(Q, A_\varepsilon(Q))$.
First, note that $d_\varepsilon- g_\varepsilon= T \circ e_\varepsilon
$, and so Lemma~\ref{lemmaexc} implies, by definition of $Q_\infty$,
that $(e_\varepsilon, d_\varepsilon- g_\varepsilon)(Q_n) \Rightarrow
(e_\varepsilon, d_\varepsilon- g_\varepsilon)(Q_\infty)$. Moreover,
$g_\varepsilon(Q_n)$ is equal in distribution to $g_\varepsilon
(\underline X_n)$ under $\PX_n^0$ and so Lemma~\ref
{lemmalevy-endpoints} shows that $g_\varepsilon(Q_n) \Rightarrow
g_\varepsilon(\underline X)$. By~\eqref{eqequality} and the definition
of $Q_\infty$, $g_\varepsilon(Q_\infty)$ and $g_\varepsilon
(\underline
X)$ are equal in distribution and so $A_\varepsilon(Q_n) \Rightarrow
A_\varepsilon(Q_\infty)$.

Let $(Q', A')$ be any accumulation point of the tight sequence $(Q_n,
A_\varepsilon(Q_n))$, so that $Q'$ is equal in distribution to $Q$ and
$A'$ to $A_\varepsilon(Q_\infty)$. Assume without loss of generality,
using Skorohod's representation theorem, that the almost sure
convergence $(Q_n, A_\varepsilon(Q_n)) \to(Q', A')$ holds: we show
that $A' = A_\varepsilon(Q')$ which will prove that $(Q', A')$ is equal
in distribution to $(Q, A_\varepsilon(Q))$.

Note $A' = (e', g', d')$: the convergence $(Q_n, A_\varepsilon(Q_n))
\to(Q', A')$ implies in view of Lemma~\ref
{lemmacontinuity-shift+stopping} and the definition of $e_\varepsilon$
that $e' = (\sigma_{d'-g'} \circ\theta_{g'})(Q')$. Since $e'$ is equal
in distribution to $e_\varepsilon(Q_\infty)$, we see that $e'$ is the
excursion of $Q'$ with endpoints $g' < d'$, and that it satisfies
$T_{e'} > \varepsilon$. To show that $e' = e_\varepsilon(Q')$ it
remains to show that this is the first such excursion. Since $Q'$ is
continuous,\vadjust{\goodbreak} it is enough to show that $\inf_{[a,a+\varepsilon]} Q' = 0$
for any $a < g'$. So let $a < g'$: since $g_\varepsilon(Q_n) \to g'$ we
must have $a < g_\varepsilon(Q_n)$ for $n$ large enough, and for those
$n$, by definition of $g_\varepsilon$~we must have $\inf_{[a,
a+\varepsilon]} Q_n = 0$. Since $\inf_{[a,a+\varepsilon]} Q_n \to
\inf_{[a,a+\varepsilon]} Q'$ by continuity, we obtain $\inf
_{[a,a+\varepsilon]} Q' = 0$ which proves that $(Q_n, A_\varepsilon
(Q_n)) \Rightarrow(Q, A_\varepsilon(Q))$. Note in particular that
since we have argued that $A_\varepsilon(Q_n) \Rightarrow
A_\varepsilon
(Q_\infty)$ we also have $A_\varepsilon(Q) = A_\varepsilon(Q_\infty)$.

\textit{Second step}. Since $Q_n$ regenerates at $0$, we
have for any measurable functions $f, h, i\dvtx D \to[0,\infty)$
\begin{eqnarray*}
&&\E_n^\varnothing\bigl( f(\sigma_{g_\varepsilon(Q_n)} Q_n)
h\bigl(e_\varepsilon(Q_n)\bigr) i(\theta_{d_\varepsilon(Q_n)}
Q_n) \bigr)
\\
&&\qquad = \E_n^\varnothing\bigl( f(\sigma_{g_\varepsilon(Q_n)}
Q_n) \bigr) \E_n^\varnothing\bigl( h
\bigl(e_\varepsilon(Q_n)\bigr) \bigr) \E_n^\varnothing
\bigl( i(Q_n) \bigr).
\end{eqnarray*}

Consider now $f, h$ and $i$ continuous and bounded, and let $n \to
+\infty$ in both sides of the previous display. Since $(Q_n,
A_\varepsilon(Q_n)) \Rightarrow(Q, A_\varepsilon(Q))$ by the first
step, Lemma~\ref{lemmacontinuity-shift+stopping} together with the
continuous mapping theorem gives
\[
\E\bigl( f(\sigma_{g_\varepsilon(Q)} Q) h\bigl(e_\varepsilon
(Q)\bigr) i(\theta
_{d_\varepsilon(Q)} Q) \bigr) = \E\bigl( f(\sigma_{g_\varepsilon
(Q)} Q) \bigr) \E\bigl(
h\bigl(e_\varepsilon(Q)\bigr) \bigr) \E\bigl( i(Q) \bigr).
\]

This implies that $\sigma_{g_\varepsilon(Q)} Q$, $e_\varepsilon(Q)$
and $\theta_{d_\varepsilon(Q)} Q$ are independent and that $\theta
_{d_\varepsilon(Q)} Q$ is equal in distribution to $Q$. Since in
addition $A_\varepsilon(Q)$ is equal in distribution to $A_\varepsilon
(Q_\infty)$ by the previous step we obtain that $\Phi_\varepsilon(Q)$
and $\Phi_\varepsilon(Q_\infty)$ are equal in distribution. For any $f
\in D$ and any $t \geq0$, one easily sees that $\Phi_\varepsilon(f)(t)
\to f(t)$ as $\varepsilon\to0$. In particular, $\Phi_\varepsilon(Q)$
converges in the sense of finite-dimensional distributions to $Q$, and
$\Phi_\varepsilon(Q_\infty)$ to $Q_\infty$, as $\varepsilon\to0$.
Hence $Q$ and $Q_\infty$ are equal in distribution which achieves the
proof of Theorem~\ref{thmconv-PS}.

\section{Identification of the limit and general initial condition}
\label{secinitial-state}

Theorem~\ref{thmconv-PS} is the most important result of the paper,
where the convergence of $(Q_n)$ under $\P_n^\varnothing$ is shown based
on the convergence of its long excursions. The formulation of
Theorem~\ref{thmconv-PS} reflects this approach, where the limiting
process is defined through its excursion measure. In general, it is not
clear whether a more explicit definition of $Q_\infty$ can be given.
For instance, in the infinite variance case we expect a similar
statement to hold, where $\Ncal$ is the excursion measure of a
reflected, spectrally positive L\'evy process; in this case we do not
know whether $Q_\infty$ can be described in another way. However, here
$\Ncal'$ turns out to be the excursion measure of $\beta^{-1}
\underline
{X}$, which is a reflected Brownian motion with drift $-\alpha/\beta$
and Gaussian coefficient $2/\beta$. This allows us to identify
$Q_\infty
$ as $\beta^{-1} \underline{X}$; see also forthcoming Theorem~\ref
{thmconv-PS-2} for a general initial condition.

%
\begin{theorem} \label{thmidentification}
$\Ncal'$ is the excursion measure of the process $\beta^{-1}
\underline X$. In particular, $Q_\infty$ is equal in distribution to
$\beta^{-1} \underline X$ under $\PX^0$.
\end{theorem}

\begin{pf}
As a variation of the original Ray--Knight theorems~\cite
{Knight630,Ray630}, it is known that $\Mcal$, the push-forward of
$\Ncal$ by the
local time process, is the excursion measure of Feller diffusion, where\vadjust{\goodbreak}
we think here of the Feller diffusion as the solution $Y$ to the
stochastic differential equation
\[
dY_t = -(\alpha/\beta) Y_t \,dt +\sqrt{(1/\beta)
Y_t} \,dB_t,
\]
with $B$ the standard Brownian motion; see, for instance, Pardoux and
Wakolbinger~\cite{Pardoux110}. It remains to show that the
push-forward of $\Mcal$ by $\Lcal$ gives $\Ncal''$, where $\Ncal''$
stands for the excursion measure of $\beta^{-1} \underline X$ (recall
that the local time of a reflected process is chosen equal to the past
infimum of the initial process).

For $\varepsilon> 0$ and $f \in D$, let $T^\varepsilon(f) = \inf\{ t
\geq0\dvtx f(t) \geq\varepsilon\}$: then $\theta_{T^\varepsilon
(\epsilon
)}(\epsilon)$, the canonical excursion shifted at time $T^\varepsilon
(\epsilon)$, under $\Ncal''( \cdot | T^\varepsilon< +\infty)$
is equal in distribution to $\sigma(\beta^{-1} X)$ under $\PX^{\beta
\varepsilon}$. The Lamperti representation theorem (see Lamperti~\cite
{Lamperti671}) asserts that $\sigma(\beta^{-1} X)$ under $\PX^{\beta
\varepsilon}$ is the image by $\Lcal$ of the Feller diffusion started
at $\varepsilon$, that is, of $\theta_{T^\varepsilon(\epsilon
)}(\epsilon
)$ under $\Mcal( \cdot | T^\varepsilon< +\infty)$. Hence the relation
\[
\Ncal'' \bigl( f \circ\theta_{T^\varepsilon} |
T^\varepsilon< +\infty\bigr) = \Mcal\bigl( f \circ\theta
_{T^\varepsilon\circ\Lcal
}
\circ\Lcal| T^\varepsilon\circ\Lcal< +\infty\bigr)
\]
holds for any nonnegative, measurable function $f\dvtx\Ecal\to
[0,\infty
)$. But by definition of $\Ncal'$ the right-hand side is precisely
$\Ncal' ( f \circ\theta_{T^\varepsilon} | T^\varepsilon<
+\infty)$ and so
\[
\Ncal'' \bigl( f \circ\theta_{T^\varepsilon} |
T^\varepsilon< +\infty\bigr) = \Ncal' \bigl( f \circ
\theta_{T^\varepsilon} | T^\varepsilon< +\infty\bigr),
\]
which can be rewritten as
\[
\Ncal'' ( f \circ\theta_{T^\varepsilon}
\indicator{T^\varepsilon< +\infty} ) = \frac{\Ncal'' (
T^\varepsilon< +\infty
)}{\Ncal' ( T^\varepsilon< +\infty)}
\Ncal' \bigl( f \circ\theta_{T^\varepsilon} \indicator
{T^\varepsilon<
+\infty} \bigr).
\]

Applying this for $f = \indicator{T^1 < +\infty}$ and $\varepsilon<
1$, we obtain
\[
\frac{\Ncal'' ( T^\varepsilon< +\infty)}{\Ncal' (
T^\varepsilon< +\infty)} = \frac{\Ncal'' ( T^1 < +\infty
)}{\Ncal' ( T^1 < +\infty)},
\]
and so for $\varepsilon< 1$ we have
\[
\Ncal'' ( f \circ\theta_{T^\varepsilon}
\indicator{T^\varepsilon< +\infty} ) = \frac{\Ncal'' ( T^1 <
+\infty)}{\Ncal
' ( T^1 < +\infty)}
\Ncal' ( f \circ\theta_{T^\varepsilon} \indicator
{T^\varepsilon<
+\infty} ).
\]

Because $f \circ\theta_{T^\varepsilon} \indicator{T^\varepsilon<
+\infty}$ converges to $f$ as $\varepsilon\to0$, we get that $\Ncal'$
and $\Ncal''$ are proportional. Moreover, since $\Ncal''$ is the
excursion measure of $\beta^{-1} \underline X$ and $\Ncal$ is that of
$\underline{X}$, we have $\Ncal''(T \in\cdot) = \Ncal(T \in\cdot)$,
and so~\eqref{eqequality} shows that the multiplicative constant must
be equal to one. This proves the result.
\end{pf}

In view of this result, it is consistent, and convenient, to redefine
$Q_\infty$ as $Q_\infty= \beta^{-1} \underline X$. In the rest of this
section, $\zeta> 0$ is some positive real number, $(\zeta_n)$ is an
integer-valued sequence such that $\zeta_n / n \to\zeta$ and we define
\[
\tau_n = \inf\bigl\{ t \geq0\dvtx L_n(0, t) >
\zeta_n / n \bigr\} \quad\mbox{and}\quad \tau= \inf\bigl\{ t \geq0\dvtx L(0,
t) = \zeta\bigr\}.
\]

The goal of this section is to prove that the sequence $(Q_n)$ under
$\P
_n^{\zeta_n*}$ converges weakly to $\beta^{-1} \underline X$ under
$\PX
^{\zeta\beta}$. To do so, we first prove that $(\sigma Q_n, T_{Q_n})$
under $\P_n^{\zeta_n^*}$ converges\vadjust{\goodbreak} weakly to $(\sigma Q_\infty,
T_{Q_\infty})$ under $\PX^{\zeta\beta}$ through a series of steps
similar to those performed in Section~\ref{secmain}. In the sequel, we
will use the fact that $S_n^*/n$ is the distribution of $X_n(\eta_n)$
under $\PX_n^0( \cdot | \eta_n < +\infty)$ where $\eta_n = \inf
\{ t \geq0\dvtx X_n(t) > 0 \}$ (see Theorem~VII.$17$ in Bertoin~\cite
{Bertoin960}); we will informally call $X_n(\eta_n)$ the overshoot of
$X_n$. The following result can be proved using standard arguments on
L\'evy processes, and so we omit the proof.

%
\begin{lemma} \label{lemmalevy-star}
The sequence $(X_n, \tau_n)$ under $\PX_n^0( \cdot | \tau_n <
+\infty)$ converges\break weakly to $(X, \tau)$ under $\PX^0( \cdot |
\tau< +\infty)$.
\end{lemma}

%
\begin{lemma}
The sequence $(Z_n)$ under $\P_n^{\zeta_n*}$ is C-tight.
\end{lemma}

\begin{pf}
By Lemma~\ref{lemmarescaling}, $Z_n = \Lcal^{-1}(Q_n)$ and so the
sequence $(Z_n)$ under $\P_n^*$ is C-tight, as can be seen by combining
Corollary~\ref{cortightness-stationary-Q} and Lemma~\ref
{lemmacontinuity+tightness-lamperti}.

Let $Z_n'$ be a process defined on the same probability space as
$Z_n$, independent of $Z_n$ and with the same law as $Z_n$ under $\P
_n^*$. Now let $Z_n''$ be the (rescaled) CMJ process defined as $Z_n'':= Z_n' + Z_n$.
Because of the lack-of-memory property of the geometric
random variable, $Z_n''$ under $\P_n^{\zeta_n*}$ has the same law as
$Z_n$ under $\P_n^*( \cdot | Z_n(0) \geq\zeta_n / n)$.

Since $(Z_n)$ under $\P_n^*$ is C-tight and $(Z_n(0))$ under $\P_n^*$
converges weakly to an exponential random variable, it is easy to show
that $(Z_n)$ under $\P_n^*( \cdot | Z_n(0) \geq\zeta_n / n)$
is C-tight. In particular, the two sequences $(Z_n'')$ and $(Z_n')$
under $\P_n^{\zeta_n*}$ are C-tight. Since the difference of two
C-tight sequences is also C-tight we obtain the C-tightness of $(Z_n)$
under $\P_n^{\zeta_n*}$ which was to be proved.
\end{pf}

%
\begin{theorem} \label{thmconv-CMJ-star}
The sequence $(Z_n)$ considered under $\P_n^{\zeta_n*}$ converges
weakly to $(L(a, \tau), a \geq0)$ under $\PX^0( \cdot | \tau<
+\infty)$. In particular, the sequence $(Z_n)$ under $\P_n^{\zeta_n*}$
converges weakly to $\Lcal^{-1}(Q_\infty)$ under~$\PX^{\beta\zeta}$,
which is the solution $Y$ to the stochastic differential equation
\[
dY_t = -(\alpha/\beta) Y_t \,dt +\sqrt{(1/\beta)
Y_t} \,dB_t,
\]
with initial condition $Y_0= \zeta$.
\end{theorem}

\begin{pf}
From Lemma~\ref{lemmarescaling} and the branching property, one gets
that $Z_n$ under $\P_n^{\zeta_n*}$ is equal in distribution to $(L_n(a,
\tau_n), a \geq0)$ considered under $\PX_n^0( \cdot | \tau_n <
+\infty)$: then the proof follows analogously as for Theorem~\ref
{thmconv-CMJ-excursion} using Lemma~\ref{lemmalevy-star} instead of
Lemma~\ref{lemmalevy-exc}. The identification of the limit as being
$\Lcal(Q_\infty)$ comes as in the proof of Theorem~\ref
{thmidentification} from a combination of the Ray--Knight theorem
together with the Lamperti representation theorem.
\end{pf}

%
\begin{lemma} \label{lemmafirst-excursion}
The sequence $(\sigma Q_n, T_{Q_n})$ under $\P_n^{\zeta_n*}$ converges
weakly to $(\sigma Q_\infty, T_{Q_\infty})$ under~$\PX^{\beta\zeta}$.
\end{lemma}

\begin{pf}
Unless otherwise specified we consider implicitly $Z_n$ and $Q_n$
under $\P_n^{\zeta_n*}$ and $Q_\infty$ under $\PX^{\beta\zeta}$.
Thanks to Corollary~\ref{corjoint}, we only have to show that $\sigma
Q_n \Rightarrow Q_\infty$ and $T_{Q_n} \Rightarrow T_{Q_\infty}$.

By Theorem~\ref{thmconv-CMJ-star}, we know that $Z_n \Rightarrow\Lcal
(Q_\infty)$. By the branching property and the fact that the overshoot
of $X_n$ is distributed like $S_n^*/n$, $T_{Z_n}$ is equal in
distribution to $\norm{X_n}_{\tau_n}$ under $\PX_n^0( \cdot |
\tau_n < +\infty)$. In view of Lemma~\ref{lemmalevy-star} combined
with the continuous mapping theorem we get weak convergence and, in
particular, tightness of $(T_{Z_n})$. Thus Lemma~\ref
{lemmacontinuity+tightness-lamperti} implies that $\Lcal(Z_n)
\Rightarrow\Lcal(\Lcal^{-1}(Q_\infty))$. Since $\Lcal(Z_n) =
\sigma
Q_n$ and $\Lcal(\Lcal^{-1}(Q_\infty)) = \sigma Q_\infty$, this proves
that $\sigma Q_n \Rightarrow\sigma Q_\infty$.

Since the workload process has the same law as $\underline X_n$,
$T_{Q_n}$ is equal in distribution to $T_{X_n}$ with the initial
condition $X_n(0) = A_n$, where $A_n$ is equal to the sum of $\zeta_n$
independent copies of $S_n^*/n$. In particular, since $\E(S_n^*) \to
\beta$ by~\eqref{eqHT}, the strong law of large number implies that
$A_n \to\zeta\beta$. Thus Lemma~\ref{lemmalevy-star} implies that
$T_{Q_n} \Rightarrow T_{Q_\infty}$ which concludes the proof.
\end{pf}

%
\begin{prop} \label{thmconv-PS-2}
The sequence $(Q_n)$ under $\P_n^{\zeta_n*}$ converges weakly to
$\beta
^{-1} \underline X$ under $\PX^{\zeta\beta}$.
\end{prop}

\begin{pf}
Let $\Ccal\dvtx D \times[0,\infty) \times D \to D$ be the concatenation
map defined for $f, h \in D$ and $s, t \geq0$ by
\[
\Ccal(f,t,h) (s) = %
\cases{ f(s), &\quad  $\mbox{if } s < t,$\vspace*{2pt}
\cr
h(t-s), & \quad $\mbox{if } s \geq t.$} %
\]

Imagine for a moment that we knew that $\Ccal$ was continuous in the
following sense: if $f_n, h_n, f, h \in D$ and $t_n, t > 0$ are such
that $f_n \to f$, $h_n \to h$, $f$ and $h$ are continuous with $f(t) =
h(t)$ and $t_n \to t$, then $\Ccal(f_n, t_n, h_n) \to\Ccal(f,t,h)$.
Then the result would follow from this result and the continuous
mapping theorem, since $Q_n = \Ccal(\sigma Q_n, T_{Q_n}, \theta Q_n)$
and $(\sigma Q_n, T_{Q_n}, \theta Q_n)$ under $\P_n^{\zeta_n*}$
converges weakly to $(\sigma Q_\infty, T_{Q_\infty}, \theta Q_\infty)$
under $\PX^{\beta\zeta}$ by Lemmas~\ref{lemmafirst-excursion} and
Theorem~\ref{thmconv-PS} [using that $(\sigma Q_n, T_{Q_n})$ and
$\theta Q_n$ are independent].

Hence we only have to prove continuity of $\Ccal$. Let $(\mu_n)$ be
any sequence of functions such that $\sup_{t \geq0}|\mu_n(t) - t|
\to
0$ and such that for each $n \geq1$, $\mu_n$ is continuous, strictly
increasing and satisfies $\mu_n(0) = 0$ and $\mu_n(t) = t_n$. Then for
any $s \geq0$ one has
\begin{eqnarray*}
&&\bigl\llvert\Ccal(f_n, t_n, h_n) \bigl(
\mu_n(s)\bigr) - \Ccal(f,t,h) (s) \bigr\rrvert\\
&&\qquad= %
\cases{
\bigl\llvert f_n\bigl(\mu_n(s)\bigr) - f(s) \bigr\rrvert,&\quad
$\mbox{if } s < t,$\vspace*{2pt}
\cr
\bigl\llvert h_n\bigl(
\mu_n(s) - t_n\bigr) - h(s-t) \bigr\rrvert,&\quad $\mbox{if } s
\geq t.$} %
\end{eqnarray*}

Since $f_n \to f$, $h_n \to h$, $f$ and $h$ are continuous and \mbox{$\sup
_{t \geq0} |\mu_n(t) - t| \to0$}, this implies that $\Ccal(f_n, t_n,
h_n) \circ\mu_n \to\Ccal(f,t,h)$. This means precisely that $(\Ccal
(f_n,t_n,h_n))$ converges in the Skorohod $J_1$ topology to $\Ccal
(f,t,h)$, and since $\Ccal(f,t,h)$ is continuous by choice of $f$ and
$h$, this means that $\Ccal(f_n, t_n, h_n) \to\Ccal(f,t,h)$. The
result is proved.
\end{pf}

\section{Discussion} \label{secdiscussion}

From a branching perspective, the result of Theorem~\ref
{thmconv-CMJ-star} is quite surprising: for the sequence $(Z_n)$ to
converge, one would naively think that the initial individuals should
start with the ``normal'' life length distribution $S_n$ instead of its
forward recurrence time $S_n^*$. This subtlety does not seem to appear
in previous works on scaling limits of continuous-time branching
processes. Although surprising from a branching perspective, this
phenomenon is well known in the folklore of queuing theory. In the rest
of this discussion, fix an integer sequence $(\zeta_n)$ such that
$\zeta
_n / n \to\zeta> 0$.

\subsection*{Discontinuity of local times}

Given some random variable $V_n$, let $Y_n^{V_n}$ be the process
obtained as follows:
\begin{itemize}
\item$Y_n^{V_n}(0)$ is distributed like $V_n/n$;
\item$Y_n^{V_n}$ has the same generator as $X_n$ in $(0,\infty)$;
\item when $Y_n^{V_n}$ hits $0$, it stays there for an exponential
duration with parameter $n \lambda_n$ and then jumps according to~$V_n/n$;
\item$Y_n^{V_n}$ is stopped at the time of its $\zeta_n$th visit to $0$.
\end{itemize}

Note also $L_n^{V_n}$ the local time process of $Y_n^{V_n}$: the
branching property together with Lemma~\ref{lemmarescaling} show that
$L_n^{V_n}$ is equal in distribution to $Z_n$ started with $\zeta_n$
individuals with i.i.d. life lengths with common distribution $V_n$.

Under the conditions that we imposed, it can be proved that $Y_n^{V_n}
\Rightarrow\underline X$ for both $V_n = S_n$ and $V_n = S_n^*$.
Nonetheless, their local time processes $(L_n^{S_n})$ and
$(L_n^{S_n^*})$ have different asymptotic behavior. On the one hand,
$(L_n^{S_n^*})$ converges in view of Theorem~\ref{thmconv-CMJ-star} to
the Feller diffusion with drift $-\alpha/\beta$ and Gaussian
coefficient $2/\beta$, started at $\zeta$. On the other hand, it can be
proved that $(L_n^{S_n})$ converges in the sense of finite-dimensional
distributions to a discontinuous process with value $\zeta$ at time~0,
but distributed for nonzero times as the same Feller diffusion started
at $\zeta/(\beta\lambda)$. In particular, the sequence $(L_n^{S_n^*})$
cannot be tight, although for each $\varepsilon> 0$ the sequence
$(\theta_\varepsilon L_n^{S_n^*})$ is tight. We now provide an
interpretation of this phenomenon in terms of state-space collapse, a
well-known property in queuing theory.

\subsection*{State space collapse}

Consider a sequence $(S_n')$ of positive random variables such that $\E
_n(S_n') \to\beta' \in(0, \infty)$, so taking for example
$S_n'=S_n^*$ yields $\beta'=\beta$ by~\eqref{eqHT}. Denote by $\P_n'$
the law of the PS queue started with $\zeta_n$ customers with i.i.d.
service requirements distributed as $S_n'$. Let $w_n$ the workload
process associated to $q_n$ and $\Wcal_n$ the rescaled process $\Wcal
_n(t) = w_n(n^2 t) / n$. Then $\Wcal_n$ is equal in distribution to
$\underline X_n$ and so converges weakly to $\underline X$, which we
write $\Wcal$ for clarity.

The state space collapse property states that the sequence $(Q_n, \Wcal
_n)$ under $\P_n^{\zeta_n *}$ converges to $(Q, \Wcal)$ which satisfy
$Q = c \Wcal$ for some constant $c > 0$. By the law of large numbers,
$\Wcal_n(0)$ under $\P_n^{\zeta_n *}$ converges to $\zeta\beta$ while
$Q_n(0)$ converges to $\zeta$, which shows that $c = \beta^{-1}$. To
understand the behavior of the processes under $\P_n'$, one needs to
zoom in around time $0$.

Define the \emph{fluid limits} $(\overline q_n)$ and $(\overline w_n)$
of $(q_n)$ and $(w_n)$ as the rescaled processes $\overline q_n(t) =
q_n(nt) / n$ and $\overline w_n(t) = w_n(nt) / n$. Fluid limits can be
thought of as functional laws of large numbers (whereas heavy traffic
approximations can be thought of as functional central limit theorems).
By definition, at the critical point the amount of work that enters the
queue is equal to the amount of work that exits it. Hence it is not
surprising that $(\overline w_n)$ under $\P_n'$ converges to the
deterministic function $w$ with constant value~$w_0 = \zeta\beta'$.
Note that the workload process does not fluctuate on the fluid time
scale $n$, while it does on the diffusion time scale $n^2$.

Let $q$ be the limit of $(\overline q_n)$ under $\P_n'$, so that $q(0)
= \zeta$; see Gromoll et~al.~\cite{Gromoll020}. Moreover it is known
that as $t$ goes to infinity, $q(t)$ converges to an equilibrium point
$q_\infty$. In steady state the residual service requirement of each
customer has mean $\beta$, which suggests, thanks to the law of large
numbers, that $q_\infty$ must satisfy $q_\infty\beta= w_0$.

So it takes a time of order of $n$ for the (scaled) queue length
process to go from $\zeta$ to $q_\infty=w_0/\beta=\zeta\beta' /
\beta
$. Since the time scale $n^2$ of the heavy traffic approximation is
orders of magnitude larger, this happens instantaneously on the
diffusion time scale and causes a discontinuity when $\beta\neq\beta'$.

Once the process has reached the equilibrium point of the fluid limit,
the state space collapse property applies. In particular, this shows
that $(Q_n)$ under $\P_n'$ should converge to a process $Q$ such that
$Q(0) = \zeta$ and $Q(t) = \beta^{-1} \Wcal(t)$ for $t > 0$. In
particular, $Q(0+) = \zeta\beta/ \beta'$ is different from $Q(0)$
when $\beta' \neq\beta$, which provides yet another interpretation of
the discontinuity of local times mentioned above. This separation of
fluid and diffusion time scales is at the heart of state space
collapse; see, for instance, Bramson~\cite{Bramson980}.

It would be interesting to understand to what extent the above
reasoning can be carried over to branching processes. In particular,
the state space collapse is quite robust in the finite variance case
and makes it possible to derive the heavy traffic limit of the PS queue
with general inter-arrival times; see Gromoll~\cite{Gromoll040}. This
suggests an approach to generalize the branching results of this paper
to the case where the offspring process is a general renewal process.

\begin{appendix}\label{app}

\section*{\texorpdfstring{Appendix: Proof of Lemma~\lowercase{\protect\ref{lemmacontinuity+tightness-lamperti}}}
{Appendix: Proof of Lemma 2.5}}

In the sequel we say that a sequence of c\`adl\`ag functions $(f_n)$ is
C-relatively compact if it is relatively compact and any of its
accumulation points is continuous. A straightforward adaptation of
Proposition~VI.$3.26$ in Jacod and Shiryaev~\cite{Jacod030}, which
gives a criterion for C-tightness, shows that a sequence $(f_n)$ is
C-relatively compact if and only if for every $m \geq0$, $\sup_{n
\geq
1} \norm{f_n}_m$ is finite and
\[
\lim_{\varepsilon\to0} \limsup_{n \to+\infty}
w_m(h_n, \varepsilon) = 0,
\]
where from now on $w_m$ is the modulus of continuity,
\[
w_m(f, \varepsilon) = \sup\bigl\{ \bigl|f(t) - f(s)\bigr|\dvtx0 \leq s,t
\leq m \mbox{ and } |t-s| \leq\varepsilon\bigr\}.
\]

%
\begin{lemma} \label{lemmaappendix-L-1}
Let $f_n, f \in\Ecal'$, and assume that the sequence $(f_n)$ is
C-relatively compact and that the sequence $(T_{f_n})$ is bounded. Then
the sequence $(\Lcal^{-1}(f_n))$ is C-relatively compact. If in
addition $f_n \to f$, then we also have $\Lcal^{-1}(f_n) \to\Lcal^{-1}(f)$.
\end{lemma}

\begin{pf}
In the rest of the proof let $h_n = \Lcal^{-1}(f_n)$ and $\overline t
= \sup_{n \geq1} T_{f_n}$. To show that the sequence $(h_n)$ is
C-relatively compact, we show that $\sup_{n \geq1} \norm{h_n}_\infty$
is finite and that
%
%
\begin{equation}
\label{eqc-rel-comp} \lim_{\varepsilon\to0} \limsup_{n \to+\infty}
w_\infty(h_n, \varepsilon) = 0
\end{equation}
with $w_\infty(j, \delta) = \lim_{m \to+\infty} w_m(j, \delta)$ for
any $j \in D$ and $\delta> 0$. By definition we have $h_n(t) =
f_n(\int_0^t h_n)$ and so $\norm{h_n}_\infty= \norm{f_n}_\infty=
\norm
{f_n}_{T_{f_n}} = \norm{f_n}_{\overline t}$. Since $(f_n)$ is
C-relatively compact, $\sup_{n \geq1} \norm{f_n}_{\overline t}$ is
finite and hence so is $\sup_{n \geq1} \norm{h_n}_\infty$.

As for $w_\infty(h_n, \varepsilon)$, we have $\int_s^t h_n \leq(t-s)
\norm{h_n}_\infty$ for any $0 \leq s \leq t$, and since $\norm
{h_n}_\infty= \norm{f_n}_{\overline t}$ we obtain for any $0 \leq s
\leq t \leq s+\varepsilon$
\[
\bigl\llvert h_n(t) - h_n(s) \bigr\rrvert= \biggl
\llvert f_n \biggl( \int_0^t
h_n \biggr) - f_n \biggl( \int_0^s
h_n \biggr) \biggr\rrvert\leq w_{\overline t} \bigl( f_n,
\varepsilon\norm{f_n}_{\overline t} \bigr).
\]
Hence $w_\infty(h_n, \varepsilon) \leq w_{\overline t} ( f_n,
\varepsilon\norm{f_n}_{\overline t} )$, and so~\eqref
{eqc-rel-comp} follows from this inequality together with the fact
that $(f_n)$ is C-relatively compact.

We now prove that $h_n \to\Lcal^{-1}(f)$ provided $f_n \to f$. Since
the sequence $(h_n)$ is C-relatively compact we only have to identify
accumulation points, so let now $h$ be any continuous accumulation
point of $(h_n)$ and assume without loss of generality that $h_n \to
h$. Let $t \geq0$: then $h_n(t) \to h(t)$, while on the other hand
from $f_n \to f$, $h_n \to h$ and the fact that $f$ is continuous we
obtain that $f_n(\int_0^t h_n) \to f(\int_0^t h)$. Since by definition
$h_n(t) = f_n(\int_0^t h_n)$, this gives $h(t) = f(\int_0^t h)$ for
every $t \geq0$. Since the solution to this equation is unique because
$f \in\Ecal'$ we obtain that $h = \Lcal^{-1}(f)$, hence the result.
\end{pf}

%
\begin{lemma} \label{lemmaappendix-L}
Let $f_n, f \in\Ecal$ and assume that $f_n \to f$, that $f$ is
continuous and that the sequence $(T_{f_n})$ is bounded. Then $\Lcal
(f_n) \to\Lcal(f)$.
\end{lemma}

\begin{pf}
Let $c_n(t) = \int_0^t f_n$ and $c_n^{-1}\dvtx[0,\infty) \to[0,T_{f_n}]$
be such that\break $\int_0^{c_n^{-1}(t)} f_n = t$ for any $t < \int f_n$ and
$c_n^{-1}(t) = T_{f_n}$ for $t \geq\int f_n$, and define similarly $c$
and $c^{-1}$ starting from $f$ instead of $f_n$. Then for any $t \geq
0$, we have by definition $\Lcal(f_n)(t) = f_n(c_n^{-1}(t))$ and
$\Lcal
(f)(t) = f(c^{-1}(t))$. We are in the framework of Theorem~$2.7$ of
Helland~\cite{Helland780}, but none of his cases applies here
(notwithstanding the problem that we allow excursions to start at $0$).
We break the proof into two steps.

\textit{First step}. Let $t < \int f$: we prove that $\norm
{c_n^{-1} - c^{-1}}_t \to0$. First, note that $c_n^{-1}$ restricted to
$[0,\int f_n]$ is the inverse of $c_n$ restricted to $[0,T_{f_n}]$, and
similarly for $c$ and $c^{-1}$. Moreover, Lemma~\ref{lemmaliminf}
implies that $T_f \leq\liminf_n T_{f_n}$ and since $\int_0^s f_n \to
\int_0^s f$ for any $s \geq0$, this implies that $\liminf_n \int f_n
\geq\int f$. Since the inverse of a continuous and strictly increasing
function is a continuous mapping (see, e.g., Theorem~$7.1$ in
Whitt~\cite{Whitt800}), we get that $\norm{c_n^{-1} - c^{-1}}_t \to0$.

\textit{Second step}. We now prove that $\Lcal(f_n) \to\Lcal
(f)$. Let $\overline t = T_f \vee\sup_n T_{f_n}$, which is finite by
assumption and is such that $c_n(t) \leq\overline t$ and $c(t) \leq
\overline t$ for any $t \geq0$ and $n \geq1$. For $t < \int f$ we write
\begin{eqnarray*}
\bigl\Vert f_n \circ c_n^{-1} - f \circ
c^{-1}\bigr\Vert_t & \leq&\bigl\Vert f_n \circ
c_n^{-1} - f \circ c_n^{-1}\bigr\Vert_t
+ \bigl\Vert f \circ c_n^{-1} - f \circ c^{-1}\bigr\Vert_t
\\
& \leq&\norm{f_n - f}_{\overline t} + w_{\overline t} \bigl(f,
\bigl\Vert c_n^{-1} - c^{-1}\bigr\Vert_t \bigr).
\end{eqnarray*}
Since $f_n \to f$, $f$ is continuous and $\norm{c_n^{-1} - c^{-1}}_{t}
\to0$, by the first step we see that the last upper bound vanishes.
Consider now some arbitrary $t' < \int f \leq t$, then
\begin{eqnarray*}
&&\bigl\Vert f_n \circ c_n^{-1} - f \circ
c^{-1}\bigr\Vert_t\\
&&\qquad \leq\bigl\Vert f_n \circ
c_n^{-1} - f \circ c^{-1}\bigr\Vert_{t'} + \sup
_{t' \leq s \leq t} f\bigl(c^{-1}(s)\bigr)
+ \sup_{t' \leq s \leq t } f_n\bigl(c_n^{-1}(s)
\bigr).
\end{eqnarray*}

The first term of this upper bound goes to $0$ by the first step. The
last term is equal to $\sup_{[c_n^{-1}(t'), \overline t]} f_n$ and
since $c_n^{-1}(t') \to c_n(t')$ by the previous step and the supremum
is a continuous function, we get
\[
\limsup_{n \to+\infty} \bigl\Vert f_n \circ c_n^{-1}
- f \circ c^{-1}\bigr\Vert_t \leq2 \sup_{c^{-1}(t') \leq s \leq T_f}
f(s).
\]

Letting $t' \to T_f$ achieves the proof.
\end{pf}

We now prove Lemma~\ref{lemmacontinuity+tightness-lamperti}, so
consider $X_n, X$ random elements of $\Ecal$ such that the sequence
$(X_n)$ is C-tight and the sequence $(T_{X_n})$ is tight. Let $(u(n))$
be any subsequence.

Assume that $X_n \Rightarrow X$: to prove that $\Lcal(X_n) \Rightarrow
\Lcal(X)$ it is enough to find a subsequence $(v(n))$ of $(u(n))$ such
that $\Lcal(X_{v(n)}) \Rightarrow\Lcal(X)$. Since the sequence
$(T_{X_n})$ is tight there exists such a subsubsequence such that
$(X_{v(n)}, T_{X_{v(n)}}) \Rightarrow(X',T')$ for some continuous $X'$
equal in distribution to $X$ and some random variable~$T'$. Lemma~\ref
{lemmaappendix-L} together with the continuous mapping theorem implies\vadjust{\goodbreak}
that $\Lcal(X_{v(n)}) \Rightarrow\Lcal(X')$, hence the result. The
other statements of the lemma follow using similar arguments, invoking
Lemma~\ref{lemmaappendix-L-1} instead of Lemma~\ref{lemmaappendix-L}.
\end{appendix}

\section*{Acknowledgments}
We would like to thank the two anonymous
referees who have contributed to significant improvements both in the
presentation and content of the paper. We are especially grateful to
one of the two referees who suggested the current proofs of Lemma~\ref
{lemmalevy-endpoints} and Theorem~\ref{thmconv-CMJ-excursion} that
lead to substantial simplifications.

%


\printaddresses

\end{document}